\documentclass[12pt]{article}
\usepackage{latexsym}
\usepackage{amssymb}
\newcommand{\rx}{\rho(x)}

\newtheorem{lemma}{Lemma}
\begin{document}

\begin{center}
{\large{\bf On the proximate order of growth of generating
functions of P\'olya frequency sequences}}\footnote{This research
has been partially supported by PRAXIS/2/2.1/FIS/286/94.}

\vspace{0.5cm}

 Maria Teresa Alzugaray

\end{center}

\abstract{We study the possible growth of $PF_r$ g.f.'s analytic
in the unit disk and describe the proximate orders of growth of
these functions. Some classical function theory results concerning
orders of growth are generalized to the case of proximate orders.}

\vspace{0.5cm}

{\bf keywords} {\it P\'olya frequency sequence, multiply positive
sequence, generating function, proximate order, total positivity}

\vspace{0.5cm}


\section{Introduction and statement of results.}

A sequence $\{a_k\}_{k=0}^\infty$ is called a P\'olya frequency
sequence of order $r$, $r\in\mathbb N\cup \{\infty\}$, also
multiply positive sequence, if all minors of order $\le r$ (all
minors if $r=\infty$) of the infinite matrix

 $$
 \left\|
  \begin{array}{ccccc}
   a_0 & a_1 & a_2 & a_3 &\ldots \\
   0   & a_0 & a_1 & a_2 &\ldots \\
   0   &  0  & a_0 & a_1 &\ldots \\
   0   &  0  &  0  & a_0 &\ldots \\
   \vdots&\vdots&\vdots&\vdots&\ddots
  \end{array}
 \right\|
 $$
are non-negative. This class of sequences is denoted by $PF_r$. We
will also denote by $PF_r$ the class of corresponding generating
functions $$ h(z)=\sum_{k=0}^\infty a_kz^k. $$ The radius of
convergence of the power series of a $PF_r$ generating function
($PF_r$ g.f.) is positive provided $r\ge2$ (\cite{tp}, p. 394).
Further we will suppose, without loss of generality, that $a_0=1$
and the radius of convergence is equal to 1, if it is not
$\infty$. We can do that because if $h(z)\in PF_r$, then $\alpha
h(\beta z)\in PF_r$ for any positive $\alpha$ and $\beta$.

It is well known (see \cite{aissen} or \cite{tp}, p. 412) that

{\bf Theorem [AESW]} {\it The class $PF_\infty$ consists of the
functions

$$ e^{\gamma z}\prod_{k=1}^\infty (1+\alpha_k z)/ (1-\beta_k z),
$$ where $\gamma \ge 0, \alpha_k\ge 0,\beta_k\ge 0$ and
$\sum(\alpha_k+\beta_k)<\infty.$}

The problem of the complete characterization of $PF_r$, $3\le
r<\infty$, has not been solved yet. The study of possible
zero-sets and growth of polynomials and entire functions belonging
to $PF_r$, $r\in\mathbb N$, has been carried out by I.J.
Schoenberg, O.M. Katkova and I.V. Ostrovskii (see \cite{schoenb},
\cite{katos} and \cite{growth}). We give a description of the
possible growth of $PF_r$ g.f.'s that are analytic in the unit
disc ${\mathbb D}$. It turns out that the growth is arbitrary in
some sense. Moreover, our construction provides a $PF_r$ g.f.
which is analytic in $\overline{\mathbb C}\backslash\{1\}$ and has
  a given growth in the  neighbourhood of $z=1$.

In this paper we will use the notion of  proximate order. Recall
that (\cite{lev}, p. 32) by definition a {\it proximate order
(p.o.)} is a function $\rho(x):\mathbb R_+\to\mathbb R_+$ which
belongs to $C^1(\mathbb R_+)$ and  satisfies the conditions:

$(i)$ $\lim\limits_{x\to\infty}\rho(x)=\rho,\ 0\le\rho<\infty;$

$(ii)$ $\lim\limits_{x\to\infty}x\frac{\rho'(x)}{\rho(x)}\log
x=0.$

We use the notation $$V(x):=x^{\rho(x)}.$$ Property $(ii)$ implies
that $V(x)$ is a strictly increasing function for  $x\ge x_0$ for
some $x_0\ge 0$. Without loss of generality, we redefine $\rx$ on
$[0,x_0)$ in such a way that $V(x)$ is strictly increasing on
$\mathbb R_+$ and $V(0)=0$.

The following properties of the functions $\rho(x)$ and $V(x)$ can
be proved by straightforward verification:\label{prop}

(a) The function $(a\rho(x)+b)/(c\rho(x)+d)\to
(a\rho+b)/(c\rho+d),\ x\to\infty,$ is a p.o. (after redefining it
on a finite interval if necessary) provided $(a\rho+b)/(c\rho+d)>
0$ and $c\rho+d\ne 0$.

(b) The function $\rho_1(x)\rho_2(x)$ is a p.o. provided
$\rho_1(x)$ and $\rho_2(x)$ are p.o..

(c) The function $\rho_1(V(x))$ is a p.o. provided $\rho_1(t)$ is
a p.o.. Moreover, if $V_1(x)=x^{\rho_1(x)}$ and
$V_2(x)=x^{\rho_2(x)}$ with p.o. $\rho_1(x)$ and $\rho_2(x)$, then
$V_1(V_2(x))=V_3(x)=x^{\rho_3(x)}$, where
$\rho_3(x)=\rho_1(V_2(x))\rho_2(x)$ is a p.o..

(d) The inverse  to $V$ function can be represented in the form
$V_{-1}(t)=t^{\tilde\rho(t)}$, where $\tilde\rho(t)\to 1/\rho,\
t\to\infty$, is a p.o., provided $\rho>0$.

(e) The function $V(x)$ is a {\it regularly varying function}
(see~\cite{sen}), that is $$
\lim_{x\to\infty}\frac{V(kx)}{V(x)}=k^{\rho}$$ uniformly on each
interval $0<a\le k\le b<\infty$. The proof of this property can be
found in \cite{lev}, p.42.

We say that the entire function $f(z)$ is of p.o. $\rho(x)$ if the
number

$$ \sigma_f=\limsup_{x\to\infty}\frac{\log M(x,f)}{V(x)} $$ is
positive and finite. Evidently, the order of the entire function
is $\rho=\lim_{x\to\infty}\rho(x)$.

Analogously, we say that the function $g(z)$ analytic in the unit
disk is of p.o. $\rx$ if the number $$
\sigma_g=\limsup_{y\to1^-}\frac{\log M(y,g)}{V(1/(1-y))} $$ is
positive and finite.

Suppose that $g$ is analytic in $\{z: 0<|z-1|\le\varepsilon\}$. We
say that $g$ has at $z=1$ a {\it singularity of p.o.} $\rho(r)$ if
the number $$ \limsup_{x\to\infty}\frac{\log M_1(x,g)}{V(x)},
$$ where
 $$M_1(x,g):=\max_{1/x\le|z-1|\le\varepsilon}|g(z)|,\;x>1/\varepsilon,$$
is positive and finite. The number $\rho=\lim_{x\to\infty}\rho(x)$
will be called the {\it order of singularity} at $z=1$.

The first main result of this paper is the following:

{\bf Theorem A} {\it Suppose that an integer $r\ge 2$ and a
proximate order $\rho(x)\to\rho,\ 0<\rho<\infty,$ are given. There
exists a function $g$ analytic in ${\mathbb D}$ and possessing the
following properties}:\\ (i) $g\in PF_r;$\\ (ii) $g$ {\it is of
p.o. $\rho(x)$ in} ${\mathbb D}$;\\ (iii) $g$ {\it is analytically
extendable to ${\overline{\mathbb C}}\backslash\{1\}$ and has at
$z=1$ an essential singularity of p.o.} $\rho(x)$.

 For the cases of $\rho=0$ and $\rho=\infty$ we
only obtained the following results:

{\bf Theorem B} {\it  For any integer $r\ge 2$ there exists a
function $\Upsilon$ analytic in $\mathbb D$ satisfying the
following properties}:\\ (i) $\Upsilon\in PF_r$;\\ (ii) $\Upsilon$
{\it is of infinite order in } $\mathbb D$;\\ (iii) $\Upsilon$
{\it is analytically extendable to ${\overline{\mathbb
C}}\backslash\{1\}$ and has at $z=1$ an essential singularity of
infinite order}.

{\bf Theorem C} {\it Suppose that an integer $r\ge 2$ and the
numbers $\rho_0,\ 1\le \rho_0 < \infty,$ and  $\sigma_0,\
0\le\sigma_0\le\infty,$ are given. There exists a function
$G(z)\in PF_r$ analytic in $\overline\mathbb C\backslash \{1\}$
and possessing an essential singularity at $z=1$ such that $$
\limsup_{y\to 1^-} \frac{\log\log M(y,G)}{\log\log
(1/(1-y))}=\rho_0
$$ and, for $\rho_0>1$, such that $$ \limsup_{y\to 1^-} \frac{\log
M(y,G)}{\left(\log (1/(1-y))\right)^{\rho_0}}=\sigma_0. $$}

The restrictions on $\rho_0$ in Theorem C are necessary. Let us
show this using the following result of \cite{dhpf}:

{\bf Theorem [DH]} {\it Let $G(z)$ be a $ PF_r$ g.f., $r\ge 2$
with radius of convergence of its power series equal to 1. Then
$(1-z)G(z)$ is a $PF_{r-1}$ g.f..}

For any function $G(z)\in PF_2$ analytic in $\mathbb D$, we have
$(1-z)G(z)\in PF_1$, therefore, $\lim_{y\to 1^-}(1-y)M(y,G)>0$ and
$\log M(y,G)>\log[1/(1-y)]+O(1),\ y\to 1^{-}.$ Thus, the quantity
$\rho_0$ in Theorem C must be greater than or equal to 1.

\section{Preliminary results and proof of Theorem~A.}

Our first lemma relates to the existence of an entire function of
a given growth belonging to $PF_r,\ r \in \mathbb N$. It is an
easy consequence of the main result of \cite{growth} due to O.M.
Katkova:

\begin{lemma}  Suppose that an integer $r\ge 2$ and a proximate
order  $\rho (x)\to\rho,\ 0\le \rho< \infty$, are given. There
exists a trascendental entire function of proximate order
$\rho(x)$ which belongs to $PF_r,\ r\in \mathbb
N$.\label{ent}\end{lemma}

The next two lemmas permit to construct $PF_r$ sequences which
will serve as sequences of Taylor coefficients of the functions
whose existence is asserted by Theorem~A.

\begin{lemma}  Suppose that $\{c_n\}_{n=0}^\infty \in PF_r,
r\ge 2, 0<\sum c_n<\infty.$ Set
$$
d_k=\sum_{n=0}^\infty\frac{c_nk^{n+r-1}}{\Gamma(n+r)},\
k=0,1,2,\cdots. $$Then $\{d_k\}_{k=0}^\infty \in PF_r$.\label{kar}
\end{lemma}

{\bf Proof} \ We derive this lemma from the following theorem due
to S. Karlin (\cite{tp}, p.107):

{\bf Theorem [Kar]} {\it Suppose that $\{c_n\}_{n=0}^\infty \in PF_r,\ r\ge 2,
\ 0<\sum c_n<\infty,\ \alpha >r-2.$ Set

$$
f_\alpha(x)=\left\{
\begin{array}{ll}
\sum_{n=0}^\infty c_n x^{n+\alpha}/\Gamma(n+\alpha+1),& x\ge 0,\\
0, x<0.
\end{array}
\right.
$$
then $f_\alpha (x)$ is a P\'olya frequency function of the order $r$.}

Recall that, $f_\alpha(x)$ is said to be a P\'olya frequency
function of order $r$, if for any $n\le r$ and for any system of
numbers $x_1<x_2<\cdots<x_n, y_1<y_2<\cdots<y_n,$ we have

$$
\det \|f_\alpha(x_j-y_i)\|_{i,j=1}^n\ge 0.
$$

Setting $\alpha=r-1,\ d_k=f_{r-1}(k), k=0,\pm 1,\pm 2, \ldots,$
and taking $x_j=j-1$ and $y_i=i-1,\ i,j\in\mathbb N$, we see that
any minor of the matrix

$$\left(
\begin{array}{ccccc}
d_0&d_1&d_2&d_3&\cdots\\
0&d_0&d_1&d_2&\cdots\\
0&0&d_0&d_1&\cdots\\
0&0&0&d_0&\cdots\\
\cdots&\cdots&\cdots&\cdots&\cdots
\end{array}
\right) $$ can be written as a minor of the matrix
$\|f_{r-1}(x_j-y_i)\|_{i,j=1}^n$. $\Box$

\begin{lemma}  Let
\begin{equation}\label{func0}f(z)=\sum_{n=0}^\infty
c_nz^n\end{equation} be an entire  $PF_r$ g.f. and set
\begin{equation}\label{fr1}
f_{r-1}(z)=\sum_{n=0}^\infty
\frac{c_n}{\Gamma(n+r)}z^{n+r-1}=\sum_{n=r-1}^\infty\frac
{c_{n-r+1}}{n!}z^n.
\end{equation} Then  $$g(z)=\sum_{k=0}^\infty
f_{r-1}(k)z^k$$ is a $PF_r$ g.f.. \label{dk}\end{lemma}

The lemma is an immediate corollary of Lemma \ref{kar}.

Our final goal is to study the growth of the function $g$ in
Lemma~\ref{dk}. For this purpose we investigate the growth of the
function $f_{r-1}(z)=:\sum_{n=0}^\infty b_nz^n$ defined by
(\ref{fr1}), where the coefficients $c_n$ are the ones in
(\ref{func0}).

For proving our next lemma we will need the theorem below, which
is due to Levin (see~\cite{lev}, p.42):

{\bf Theorem [Levin]} {\it Let $\rho(x)\to\rho>0,\ x\to\infty$, be
a p.o.. For the entire function $B(z)=\sum_{n=0}^\infty b_nz^n$
set $$ \sigma_B=\limsup_{n\to\infty}\frac{\log M(x,B)}{V(x)}. $$
Then $$
(\sigma_Be\rho)^{1/\rho}=\limsup_{n\to\infty}V_{-1}(n)\sqrt[n]{|b_n|}.
$$}

\begin{lemma}\label{ord}
If the function (\ref{func0}) is an entire function of p.o.
$\rho(x)$, then the function (\ref{fr1}) is an entire function of
p.o.
\begin{equation}\label{rho1}
\rho_1(x)=\frac{\log \psi(x)}{\log x}\to\rho_1=\frac\rho{\rho+1},\
x\to\infty,
\end{equation}
where $\psi(x)$ is the inverse function of $tV_{-1}(t)$ and
$V_{-1}(t)$ is the inverse function of $V$.\end{lemma}

{\bf Proof} By (\ref{rho1}) we have $\psi(x)=x^{\rho_1(x)}$, and
$(\psi(x))_{-1}=tV_{-1}(t)$. By the properties of p.o.'s d) and a)
listed on p. \pageref{prop}, we have that $\rho_1(x)$ is a p.o.

For proving that $f_{r-1}(z)$ is of p.o. $\rho_1(x)$ we calculate
$$ \limsup_{n\to\infty}\psi_{-1}(n)\sqrt[n]{|b_n|}, $$where $b_n$
are the Taylor coefficients of $f_{r-1}$.

We have $$ \limsup_{n\to\infty}\psi_{-1}(n)\sqrt[n]{|b_n|}=
\limsup_{n\to\infty}nV_{-1}(n)\sqrt[n]{c_{n-r+1}/n!}= $$
$$=e\limsup_{n\to\infty}V_{-1}(n)\sqrt[n]{c_{n-r+1}}=
e\limsup_{n\to\infty}V_{-1}(n)\sqrt[n]{c_{n}}.$$

By Levin's theorem $$
\limsup_{n\to\infty}\psi_{-1}(n)\sqrt[n]{|b_n|}=e(\sigma_f
e\rho)^{1/\rho} $$ is positive and finite, thus $\rho_1(x)$ is the
p.o. of $f_{r-1}(z)$. $\Box$

In the sequel, we will need a number of facts concerning functions
\begin{equation}
h(z)=\sum_{k=0}^\infty a_kz^k
\label{func}
\end{equation}
which are analytic in $\mathbb D$.

The following lemma is an analogue of Levin's theorem for the case
of functions analytic in the unit disk.

\begin{lemma}\label{levin2} Let $\rho(x)$ be a proximate order, $\rho(x)\to\rho>0,\
x\to\infty$, and $\xi(t)$ the inverse function of $xV(x)$. For the
function (\ref{func}) set $$ \sigma_h=\limsup_{y\to 1^-}\frac
{\log M(y,h)}{V(1/(1-y))}. $$ Then $$
\frac{\rho+1}{\rho}(\sigma_h\rho)^{1/(\rho+1)}=
\limsup_{k\to\infty}\frac{\xi(k)}{k}\log|a_k|. $$ \end{lemma}

We could not find in the literature neither this result, nor the
Lemmas \ref{x} and \ref{xx}, which will appear later, in spite of
the fact that there are plenty of similar results with similar
proofs (see, for example \cite{she1}, \cite{she2} where one can
find further bibliography). For the reader's convenience we
present the proofs in the last section of the paper.

{\bf Theorem [Wigert]} (see, for example \cite{lev}, p. 394) {\it
The function $h(z)$ defined by $(\ref{func})$ can be analytically
extended to $\overline\mathbb C\backslash\{1\}$ and is equal to
zero at infinity if and only if there exists an entire function
$A(z)$ whose growth is not greater than of order 1 and minimal
type such that $a_k=A(k), k=0,1,2,\ldots.$}

As an improvement to Wigert's theorem Faber established the
relation between the orders of growth of the functions $A(z)$ and
$h(z)$ (see~\cite{faber} or \cite{bieb}, \S 1, Th. 1.3.11):

{\bf Theorem [Faber]} {\it The function $h(z)$ can be analytically
extended to $\overline\mathbb C\backslash\{1\}$  and has at $z=1$
a singularity of order $\rho$ if and only if the entire function
$A(z)$ is of order of growth $\rho_A=\rho/(\rho+1)$.}

We were able to relate the p.o. of growth of the function $h$ in
$\mathbb D$ with the p.o. of growth of the function $A$ when $A$
has non-negative  Taylor coefficients.

\begin{lemma}\label{gelf2} Let $A(z)$ be an entire function of proximate order
$\rho_1(x)\to\rho_1,\ 0<\rho_1<1,\ x\to\infty$, with non-negative
Taylor coefficients, such that $A(n)=a_n,\ n=0,1,2,\ldots$, where
$a_n$ are the coefficients of the function (\ref{func}). Then
$h(z)$ is of p. o.
\begin{equation}\label{rogelf2}\rho(x)=\frac{\rho_1(\xi_{-1}(x))}
{1-\rho_1(\xi_{-1}(x))}\end{equation} in $\mathbb D$, where
\begin{equation}\label{xigelf2}\xi(t)=t^{1-\rho_1(t)}
\end{equation}
is the inverse of the function $xV(x)$. Moreover, $$
\rho_1(t)=\frac{\rho(\xi(t))}{\rho(\xi(t))+1}.$$ Also, $A(z)$ is
of order 0 if and only if $h(z)$ is of order 0 in $\mathbb D$.
\end{lemma}

{\bf Proof} \ By properties (a), (d) and (c) of proximate orders
(see p.~\pageref{prop}), we can conclude that $\rho(x)$ is a p.o..

Now we want to prove that $\rho(x)$ is the p.o. of (\ref{func}),
i.e. that the number $$ \sigma_h=\limsup_{y\to 1^-}\frac{\log
M(y,h)}{V(1/(1-y))} $$ is positive and finite. We will prove it
with the help of Lemma \ref{levin2}. Let $\xi(x)$ and $\rho(x)$ be
defined by (\ref{xigelf2}) and (\ref{rogelf2}), respectively. We
have $\rho(\xi(t))=\rho_1(t)/(1-\rho_1(t))$.  Therefore,
$t=\xi(t)V(\xi(t))$, which means that $\xi(t)$ is the inverse
function of $xV(x)$.

On the other hand, denoting $V_1(x)=x^{\rho_1(x)}$, we get by
(\ref{xigelf2}) $\xi(n)/n=1/V_1(n)$ and, hence, $$
\limsup_{n\to\infty}\frac{\xi(n)}{n}\log a_n=
\limsup_{n\to\infty}\frac{\log A(n)}{V_1(n)}=
\limsup_{x\to\infty}\frac{\log M(x,A)}{V_1(x)}, $$ which is
positive and finite. Note that we were able to write the last
equality because $A(z)$ has non-negative Taylor coefficients. Then
by Lemma \ref{levin2}, $\rho(x)$ is the proximate order of $h(z)$.

The last assertion of the lemma follows from the following result
of Beuermann (see~\cite{beuer}):

{\bf Theorem [Beuermann]}{\it For the function $h(z)$, set
$$
\lambda=\limsup_{y\to 1^-}\frac {\log\log M(y,h)}{\log (1/(1-y))}.
$$
Then the following equality is valid}
$$
\frac\lambda{\lambda+1}
=\limsup_{k\to\infty}\frac{\log^+\log^+|a_k|}{\log k }.
$$
Thus, by Beuermann's theorem $A(z)$ being of order 0 is equivalent
to
$$0=\limsup_{x\to\infty}\frac{\log\log M(x,A)}{\log x }
=\limsup_{k\to\infty}\frac{\log\log A(k)}{\log k}$$ $$
=\limsup_{k\to\infty}\frac{\log^+\log^+|a_k|}{\log k }
=\limsup_{y\to 1-}\frac{\log\log M(y,h)}{\log (1/(1-y)) }.$$
 $\Box$

The lemma below is a generalization of the theorem of Faber to the
case of proximate orders that do not tend to zero at infinity.

\begin{lemma}\label{faber2} Let $h(z)$ be analytically extendable
to $\overline\mathbb C\backslash\{1\}$ with a singularity at $z=1$
of p.o. $\rho(x)\to\rho>0,\ x\to\infty$, and $A(w)$ be of p.o.
$\rho_A(x)\to\rho_A<1,\ x\to\infty$. Then
\begin{equation}\label{rhoA}
\rho_A(t)=\frac{\rho(\xi(t))}{\rho(\xi(t))+1}, \end{equation}
where $\xi(t)$ is the inverse function of $xV(x)$. Moreover,

\begin{equation}\label{rofaber2}\rho(x)=\frac{\rho_A(\xi_{-1}(x))}
{1-\rho_A(\xi_{-1}(x))}\ {\rm and}\
\xi(t)=t^{1-\rho_A(t)}.\end{equation}
\end{lemma}

{\bf Proof} Our proof is based on the method used by Gelfond (see
\cite{gelf}) to obtain Faber's result.

By hypothesis $h(z)$ has at $z=1$ an essential singularity of p.o.
$\rho(x)$. Thus,  there exists a $\sigma$ and a $x_\sigma$ such
that
\begin{equation}\label{maxh}
\max\left\{|h(z)|:|1-z|=\frac 1x\right\}<\exp\{\sigma V(x)\}
\end{equation}
for $x>x_\sigma$.

On the other hand, since $A(k)=a_k,\ k=0,1,2,\ldots,$ we can write
(see, for example, \cite{bieb}, \S 1, equation (1.3.29))
\begin{equation}\label{integr}
A(w)=-\frac 1{2\pi i}\int\limits_{|z-1|=1/x}h(z)\exp\{-(w+1)\log
z\}dz,
\end{equation}
where we take $x>x_\sigma$. Consider also $x\ge 2$, so $|\arg
z|\le 2/x$. We have
\begin{equation}\label{expo}
|\exp\{-(w+1)\log z\}|<\exp\left\{|w+1|\left(|\log|z||+\frac
2x\right)\right\}\le \exp\left\{\frac {5|w|}x\right\}
\end{equation}
for $|1-z|=1/x$ and $|w|$ large enough.

Denoting $|w|$ by $t$, and using (\ref{integr}), (\ref{maxh}) and
(\ref{expo}) we can write $$ |A(w)|\le \exp\left\{\sigma
V(x)+\frac{5t}x\right\},$$ for $t$ and $x$ large enough.

In the previous inequality we set $x:=\xi(5t/(\rho\sigma))$, where
$\xi(t)$ is the inverse function of $xV(x)$. Note that we can do
that since $\xi(t)$ is an increasing function (see properties of
p.o.'s on p. \pageref{prop}). We obtain $$ |A(w)|\le
\exp\left\{\sigma
V(\xi({5t}/({\rho\sigma})))+\frac{5t}{\xi({5t}/({\rho\sigma}))}\right\}$$
$$=\exp\left\{5\left(\frac{\rho+1}\rho\right)\frac{t}{\xi({5t}/({\rho\sigma}))}\right\}
$$
$$=\exp\left\{5\left(\frac{\rho+1}\rho\right)\frac{\xi(t)V(\xi(t))}{\xi({5t}/({\rho\sigma}))}
\right\}.$$

But, since $\xi(t)$ is a regularly varying function, $$
\xi(\frac{5t}{\rho\sigma})=\left(\frac
5{\rho\sigma}\right)^{1/(\rho+1)}\xi(t)\{1+o(1)\},\ t\to\infty. $$

Thus, asymptotically $$
|A(w)|\le\exp\left\{K\sigma^{1/(\rho+1)}V(\xi(t))\{1+o(1)\}\right\},\
t\to\infty, $$ where $K=K(\rho)>0$ is a constant.

Therefore,
\begin{equation}\label{maj1}
\limsup_{t\to\infty}\frac{\log
M(t,A)}{V(\xi(t))}\le(\sigma_H)^{1/(\rho+1)}<\infty.
\end{equation}

Conversely, suppose that $\rho_A(t)$ is the p.o. of $A(w)$, i.e.
$$\sigma_A=\limsup_{t\to\infty}\frac{\log M(t,A)}{V_A(t)},$$ where
$V_A(t)=t^{\rho_A(t)}$, is positive and finite. Then for
$\sigma>\sigma_A$ the following asymptotic inequality holds
\begin{equation}\label{vat}
|A(w)|<\exp\{\sigma V_A(t)\}.
\end{equation}

On the other hand, we have $h(z)=H(1/(1-z))$, where $H$ is an
entire function. Moreover, $$
M_1(x,h)=\max\{|h(z)|:\varepsilon>|z-1|\ge 1/x\}$$
$$=\max\{|H(w)|: |w|\le x\}=M(x,H),$$ so the p.o.'s of the function
$H$ and of the essential singularity $z=1$ of $h$ coincide. Also,
$h(z)=\tilde H(z/(1-z))$ and the inequality
$$ M(t-1,\tilde H)\le M(t,H)\le M(t+1,\tilde H)$$ shows that $H$
and $\tilde H$ are of the same p.o..

Denoting $z/(1-z)$ by $\zeta$ we have that  $$
h(z)=\sum_{k=0}^\infty a_kz^k=(1+\zeta)\sum_{n=0}^\infty
h_n\zeta^n.$$

Note that the p.o.'s of $H(w)$ and $\sum_{n=0}^\infty h_n\zeta^n$
coincide. Next we calculate the coefficients $h_n$: $$
h(z)=(1+\zeta)\sum_{k=0}^\infty a_k\frac{\zeta^k}{(1+\zeta)^{k+1}}
$$ $$=(1+\zeta)\sum_{k=0}^\infty a_k\zeta^k\left(\sum_{m=0}^\infty
{m+k\choose k}(-1)^m\zeta^m\right)$$
$$=(1+\zeta)\sum_{n=0}^\infty\left(\sum_{k=0}^n
(-1)^{n-k}{n\choose k}a_k\right)\zeta^n$$

Thus, we have $$h_n=\sum_{k=0}^n(-1)^{n-k}{n\choose k}a_k.$$

On the other hand, it is not difficult to see that
$$\frac{A(w)}{w(w-1)\ldots(w-n)}=\sum_{k=0}^n\frac{\alpha_k}{w-k},$$
where
$$\alpha_k=\left.\frac{(w-k)A(w)}{w(w-1)\ldots(w-n)}\right|_{w=k}=\frac{1}{n!}{n\choose
k}(-1)^{n-k}a_k.$$

Therefore, by the residue theorem we can write
$$h_n=\frac{n!}{2\pi i}
\int\limits_{|w|=R}\frac{A(w)}{w(w-1)...(w-n)}dw,$$ for $R>n$.

It is obvious that $$ |h_n|\le
n!\frac{M(R,A)}{(R-1)...(R-n)}=n!\frac{\Gamma(R-n)}{\Gamma(R)}M(R,A),\
R>n.$$

Taking $R$ large enough in the previous inequality, we obtain by
(\ref{vat}) that $$ |h_n|\le
n!\frac{\Gamma(R-n)}{\Gamma(R)}\exp\{\sigma V_A(R)\}.$$

Thus,
\begin{equation}\label{gamma}
 |h_n|^{1/n}\le
\left[\frac{n!\Gamma(R-n)}{\Gamma(R)}\right]^{1/n}\exp\{\frac\sigma
n V_A(R)\}.
\end{equation}

Now we denote by $R(t)$ the function inverse to $V_A(x)$ and set
$R:=R(n)$ on the right-hand side of (\ref{gamma}). We can do this
because
 $R(t)=t^{\alpha(t)},$ where $\alpha(t)$ is
a p.o. such that $\alpha(t)\to 1/\rho_A>1,\ t\to\infty$ (see
properties on p. \pageref{prop}). So
$\lim_{n\to\infty}R(n)/n=\infty$.

We apply the Stirling formula to the right hand of (\ref{gamma}).
Note that $$ (R(n))^{1/n}=n^{\alpha(n)/n}\sim 1,\ n\to\infty,$$
so, when $n\to\infty$,
$$[\Gamma(R)]^{1/n}\sim
\left(\frac{R(n)}{e}\right)^{R(n)/n}$$

\begin{equation}\label{stir1}
=\left(\frac{n^\alpha}{e}\right)^{n^{\alpha-1}}=\exp\{\alpha
n^{\alpha-1}\log n-n^{\alpha-1}\},
\end{equation}
where $\alpha=\alpha(n)$.

Also, when $n\to\infty,$ $$[\Gamma(R-n)]^{1/n}\sim
\left(\frac{R(n)-n}{e}\right)^{R(n)/n-1}$$
$$=\left(\frac{n^\alpha-n}{e}\right)^{n^{\alpha-1}-1}=
\left(\frac{n^\alpha}{e}\right)^{n^{\alpha-1}-1}(1-n^{1-\alpha})^{n^{\alpha-1}-1}
$$ $$=\exp\{\alpha(n^{\alpha-1}-1)\log
n-n^{\alpha-1}+1\}(1-n^{1-\alpha})^{n^{\alpha-1}-1}$$
\begin{equation}\label{stir2}
\sim\exp\{\alpha(n^{\alpha-1}-1)\log n-n^{\alpha-1}\},
\end{equation}
where $\alpha=\alpha(n)\to 1/\rho_A>1,\ n\to\infty.$

 By (\ref{stir1}) and (\ref{stir2}), we have that the right
hand of (\ref{gamma}), when $n\to\infty,$ is equivalent to $$\frac
ne\exp\{\alpha(n^{\alpha-1}-1)\log n-n^{\alpha-1}-\alpha
n^{\alpha-1}\log n+n^{\alpha-1}+\sigma\}$$ $$= e^{\sigma-1}\frac
n{R(n)}\{1+o(1)\},\ n\to\infty.$$

Thus, $$\limsup_{n\to\infty}\frac{R(n)}n|h_n|^{1/n}<\infty.$$ Note
that $\log[(R(t)/t)_{-1}]/\log x,$ where $(R(t)/t)_{-1}(x)$ is the
inverse function to $R(t)/t$, is a p.o.. By Levin's Theorem we
have
\begin{equation}\label{maj2}
\limsup_{x\to\infty}\frac{\log M(x,\tilde
H)}{(R(t)/t)_{-1}(x)}<\infty.
\end{equation}

We will say that $\tilde H$ is "not greater than"
$(R(t)/t)_{-1}(x).$

Now we will prove that the quantities we refer to in (\ref{maj1})
and (\ref{maj2}), that is $\limsup_{t\to\infty}{\log
M(t,A)}/{V(\xi(t))}$ and $\limsup_{x\to\infty}{\log M(x,\tilde
H)}/{(R(t)/t)_{-1}(x)}$, are also both positive.

For $V_A(t)$ we can find a $V(x)=x^{\rho(x)}$ such that
$V_A(t)=V(\xi(t))$ and $\xi(t)=(xV(x))_{-1}$. Indeed, setting
$\xi(t)=t/V_A(t)$ and $V(x)=(\xi(t))_{-1}/x$, we have
$\xi(t)=(xV(x))_{-1}$ and $V_A(t)=t/\xi(t)=V(\xi(t))$. We will
prove that $\rho(x)$ is the p.o. of $h(z)$.

By what was proved earlier, $\tilde H$ is "not greater than"
$(R(t)/t))_{-1}$, $ R(t)=(V_A(x))_{-1}.$ But
$R(t)=(V_A(x))_{-1}=(V(\xi(x)))_{-1}=\xi_{-1}(V_{-1}(t))=tV_{-1}(t)$
because $\xi(t)$ is the inverse function to $xV(x)$. Thus, $$
\limsup_{t\to\infty}\frac {\log M(t,\tilde H)}{V(t)}<\infty.$$

Now suppose that  $$\limsup_{t\to\infty}\frac {\log M(t,\tilde
H)}{V(t)}=0,$$ then, by (\ref{maj1}), $$\limsup_{t\to\infty}\frac
{\log M(t,A)}{V(\xi(t))}=0,$$ which contradicts the fact
$V(\xi(t))=V_A(t)$.

We have necessarily $$0<\limsup_{t\to\infty}\frac {\log M(t,\tilde
H)}{V(t)}<\infty$$  and $h(z)$ is a function of p.o. $\rho(x)=\log
V(x)/\log x$.

Note that $$\rho_A(t)=\frac{\log V(\xi(t))}{\log
t}=\rho(\xi(t))\frac{\log\xi(t)}{\log
t}=\frac{\rho(\xi(t))}{\rho(\xi(t))+1}$$ and hence
$$\rho(x)=\frac{\rho_A(\xi_{-1}(x))}{1-\rho_A(\xi_{-1}(x))}.$$ By
construction $\xi(t)=(xV(x))_{-1}=t/V_A(t)$. $\Box$

The following result is closely related to the main result of
\cite{macint} due to Macintyre and Wilson:

\begin{lemma}\label{igualdad}
Let $h$ be a function analytic in $\mathbb D$ of order greater
than 0 whose Taylor coefficients are interpolated by the entire
function $A$ of order less than 1 and having non-negative Taylor
coefficients. Then the p.o. of $h$ in $\mathbb D$ and the p.o. of
its singularity at $z=1$ coincide. Also, $h$ is of order 0 in
$\mathbb D$ if and only if its singularity at $z=1$ is of order 0.
\end{lemma}

The lemma is an immediate consequence of Lemmas \ref{gelf2} and
\ref{faber2} and Faber's theorem.

{\bf Proof of Theorem A} Let $r\ge 2$ and a p.o. $\rho(x)$ be
given. Let $f(z)$ be the entire function whose existence is
established by Lemma~\ref{ent}. Set
$$ g(z)=\sum_{k=0}^\infty f_{r-1}(k)z^k,
$$where $f_{r-1}$ is defined by (\ref{fr1}). According to
Lemma \ref{dk} we have $g(z)\in PF_r$.

By Lemma \ref{ord} the function $f_{r-1}(z)$ is an entire function
of order $\rho_1=\rho/(\rho+1)<1$. According to Lemma \ref{gelf2}
$g(z)$ is of p.o.
\begin{equation}\label{portilde}
\tilde\rho(x)=\frac{\rho_1(\xi_{-1}(x))}{1-\rho_1(\xi_{-1}(x))}
\end{equation}
in $\mathbb D,$ where $\rho_1(t)$ is the p.o. of the function
$f_{r-1}(z)$ and
\begin{equation}\label{eq*}\xi(t)=t^{1-\rho_1(t)}.
\end{equation}

Now we will show that $\tilde\rho(x)=\rho(x)$.

By virtue of Lemma \ref{ord} we have

\begin{equation}\label{eqpsi}
\rho_1(x)=\frac{\log\psi(x)}{\log x},
\end{equation}
where $\psi(x)$ is the inverse function of $tV_{-1}(t),\
V(x)=x^{\rho(x)},$ and $\rho(x)$ is the p.o. of the function
$f(z)$.

For the function $\xi(t)$ defined by (\ref{rofaber2}) the
following equality holds
\begin{equation}\label{int}
\xi(\psi_{-1}(t))=V_{-1}(t).
\end{equation}

Indeed, by (\ref{eq*}), (\ref{eqpsi}) and the definition of $\psi$
we can write

$$
\xi(\psi_{-1}(t))=(\psi_{-1}(t))^{1-\rho_1(\psi_{-1}(t))}=\psi_{-1}(t)^{1-\log
t/\log\psi_{-1}(t)}= $$ $$=\exp\{(1-\frac{\log
t}{\log\psi_{-1}(t)})\log\psi_{-1}(t)\}=\exp\{\log\psi_{-1}(t)-\log
t\}=V_{-1}(t). $$

Now with the aid of (\ref{int}), (\ref{portilde}), (\ref{eqpsi})
and the definition of $V$ we are able to write

$$ \tilde\rho(V_{-1}(t))=\tilde\rho(\xi(\psi_{-1}(t)))=
\frac{\rho_1(\psi_{-1}(t))}{1-\rho_1(\psi_{-1}(t))}= $$ $$
=\frac{\log t}{\log\psi_{-1}(t)-\log t}=\frac{\log t}{\log
V_{-1}(t)}=\frac{\log V(V_{-1}(t))}{\log V_{-1}(t)}=\rho
(V_{-1}(t)). $$

Thus, we have shown that $\tilde\rho(x)=\rho(x)$.

On the other hand, the Wigert theorem is applicable to the
function $g$ and it shows that $g(z)$ can be extended to
$\overline{\mathbb C}\backslash\{1\}$. According to Lemma
\ref{igualdad} $g(z)$ has at $z=1$ an essential singularity of
p.o. $\rho(x)$.$\Box$

\section{Preliminary results and proof of Theorem B.}

The following result due to O.M. Katkova and I.V. Ostrovskii will
be useful to prove Theorem B.

{\bf Theorem [KatOst]} {\it Let $g_1$ be an arbitrary entire
function which is positive on the positive $x$-axis and such that
$g_1(0)=1$. For every $r\in\mathbb N$ there exists an entire
function $g_2\in PF_r$ such that $g_2g_1\in PF_r$.}

\begin{lemma}\label{infinito}
There exists an entire function of infinite order belonging to
$PF_r,\ r\in\mathbb N.$
\end{lemma}

This lemma is a consequence of Theorem [KatOst]. Indeed, let
$g_1=\sum_{n=0}^\infty \varphi_nz^n$ be an entire function of
infinite order with positive coefficients. By Theorem [KatOst]
there exists an entire function $g_2$ such that $\Phi=g_2g_1\in
PF_r$. Evidently, $\Phi(z)$ is of infinite order.

{\bf Proof of Theorem B} Let the integer $r\ge 2$ be given. Let
$$
\Phi(z)=\sum_{n=0}^\infty\phi_n z^n $$ be the entire $PF_r$ g.f.
of infinite order whose existence is established by Lemma
\ref{infinito}. Setting $c_n=\phi_n$, $f(z)=\Phi(z)$ and
$$f_{r-1}(z)=\Phi_{r-1}(z)=\sum_{n=0}^\infty
\frac{\phi_nz^{n+r-1}}{(n+r)!}$$ in Lemma 3, we obtain that the
function
$$\Upsilon(z)=\sum_{k=0}^\infty \Phi_{r-1}(k)z^k $$ is a $PF_r$ g.f..

Using Hadamard and Stirling formulas we will show that the entire
function $\Phi_{r-1}(z)$ is of order 1. Indeed, its order is equal
to $$ \limsup_{n\to\infty}\frac{n\log n}{\log(n!/\phi_{n-r+1})}=
\left[\liminf_{n\to\infty}\left(\frac{\log(n!)}{n\log
n}+\frac{\log(1/\phi_{n-r+1})}{n\log n}\right)\right]^{-1}=1,$$
since $\Phi(z)$ is of infinite order, which means that
$$\limsup_{n\to\infty}\frac{n\log
n}{\log(1/\phi_{n-r+1})}=\infty.$$ Moreover, $\Phi_{r-1}(z)$ is of
minimal type. Indeed, for any $\varepsilon>0$ we have the
asymptotic inequality $\phi_n<\varepsilon^n$ and, thus,
$M(x,\Phi_{r-1})=\Phi_{r-1}(x)\le O(x^{r-1}\exp\{\varepsilon
x\})$, for $x\to\infty$.

Suppose that  $\Upsilon(z)$ is of finite order $\rho$ in $\mathbb
D$, then by Lemma \ref{gelf2} the function $\Phi_{r-1}$ is of
order of growth equal to $\rho/(\rho+1)<1$ which is a
contradiction. Hence, the order of growth of $\Upsilon$ in
$\mathbb D$ is infinite.

The Wigert theorem can be applied to $\Upsilon(z)$ showing that
this function can be analytically continued to $\overline\mathbb
C\backslash\{1\}$. By Lemma \ref{igualdad} $\Upsilon(z)$ has an
essential singularity of infinite order at $z=1$. $\Box$

\section{Preliminary results and proof of Theorem C.}

For proving Theorem C we will need the lemmas below which are
analogs of Hadamard formulas connecting the growth of a function
with its Taylor coefficients. We present the proofs of these
lemmas in the last section of the paper.

\begin{lemma}  Let $F(z)=\sum_{n=0}^\infty C_nz^n$ be an entire
trascendental function. Set
$$
 \rho_0=\limsup_{x\to
\infty}\frac {\log\log M(x,F)}{\log\log x} $$
 and, for $1<\rho_0<\infty$, set
 $$ \sigma_0=\limsup_{x\to
\infty}\frac {\log M(x,F)}{(\log x)^{\rho_0}}. $$
 Then $$
\frac{\rho_0^{\rho_0}\sigma_0}{(\rho_0-1)^{\rho_0-1}}
=\limsup_{n\to\infty}\frac{n^{\rho_0}}{\left(\log(1/|C_n|)\right)^{\rho_0-1}}.$$
\label{x}\end{lemma}

We say that the function $F(z)$ of Lemma \ref{x} is of {\it
logarithmic order} $\rho_0$ and {\it logarithmic type} $\sigma_0$.

\begin{lemma} Let $h(z)=\sum_{k=0}^\infty a_kz^k$ be a function
analytic in $\mathbb D$. Set $$ \rho_0=\limsup_{y\to 1^-}\frac
{\log\log M(y,h)}{\log\log (1/(1-y))}. $$ Let $\rho_0$ satisfy
$1\le\rho_0<\infty$, and for $\rho_0>1$ set
 $$ \sigma_0=\limsup_{y\to
1^-}\frac {\log M(y,h)}{\left(\log (1/(1-y))\right)^{\rho_0}}. $$
 Then $$ \rho_0=
\limsup_{k\to\infty}\frac{\log^+\log^+|a_k|}{\log\log k} {\rm\ \ \
and\ \ \ } \sigma_0=\limsup_{k\to\infty}\frac{\log^+|a_k|}{(\log
k)^{\rho_0}}.$$\label{xx}\end{lemma}

{\bf Proof of Theorem C} Applying Lemma \ref{ent} with
$\rho(x)=\rho_0\log\log x/\log x,$  $V(x)=x^{\rho(x)}$ we have
that for any given $r\in\mathbb N$ there exists an entire
trascendental function$$ F(z)=\sum_{n=0}^\infty C_nz^n $$ of
logarithmic order $\rho_0$ and of logarithmic type $\sigma_0$
belonging to $PF_r$. Setting $c_n=C_n$, $f(z)=F(z)$ and
$f_{r-1}(z)=F_{r-1}(z)=\sum_{n=0}^\infty C_nz^{n+r-1}/\Gamma(n+r)$
in Lemma \ref{dk} we obtain that the function
 $$ G(z)=\sum_{k=0}^\infty F_{r-1}(k)z^k=\sum_{k=0}^\infty D_kz^k $$
is a $PF_r$ g. f..

Note that the logarithmic order and type of the entire functions
$F(z)$ and $F_{r-1}(z)$ coincide. Indeed, $$ \rho_0=\limsup_{x\to
\infty}\frac{\log\log M(x,F)}{\log\log x}. $$ Therefore, by Lemma
\ref{x}, we have $$ \frac{\rho_0 -1}{\rho_0}=
\limsup_{n\to\infty}\frac{\log n}{\log\log (1/C_n)}. $$ Noting
that $$ \liminf_{n\to\infty}\frac {\log
(1/{C_n})}{n^{1+\alpha}}\ge 1,\ {\rm for\ some}\ \alpha>0, $$ we
can calculate $$ \limsup_{n\to\infty}\frac{\log n}{\log\log
(\Gamma(n+r)/C_n)}=\frac {\rho_0-1}{\rho_0}. $$ Hence, applying
Lemma \ref{x}, $$ \rho_0=\limsup_{x\to \infty}\frac {\log\log
M(x,F_{r-1})}{\log\log x}. $$ Analogously, the logarithmic types
coincide because $$\limsup_{n\to\infty}\frac{n^{\rho_0}}{(\log
(\Gamma(n+r)/C_n))^{\rho_0-1}}=\limsup_{n\to\infty}\frac{n^{\rho_0}}{(\log
(1/C_n))^{\rho_0-1}}.$$

Hence, the function $F_{r-1}(z)$ is an entire function of order 0.
By Wigert's Theorem $G(z)$ can be extended to $\overline\mathbb
C\backslash\{1\}$.

Since the coefficients of $F_{r-1}(z)$ are non-negative we have
$F_{r-1}(x)=M(x,F_{r-1})$ for $x>0$ and therefore, remembering
that $D_k=F_{r-1}(k)$ we can write
 $$
\rho_0=\limsup_{t\to\infty}\frac{\log\log F_{r-1}(t)}{\log\log t
}=\limsup_{k\to\infty}\frac{\log \log D_k}{\log\log k } $$ and $$
\sigma_0=\limsup_{t\to\infty}\frac{\log F_{r-1}(t)}{(\log
t)^{\rho_0} }=\limsup_{k\to\infty}\frac{\log D_k}{(\log
k)^{\rho_0}},\ {\rm for} \ \rho_0>1.$$
 Applying Lemma \ref{xx} to $G(z)$ we can assert that
$$ \limsup_{y\to 1^-}
\frac{\log\log M(y,G)}{\log\log (1/(1-y))}=\rho_0$$ and
 $$ \limsup_{y\to
1^-}\frac {\log M(y,h)}{\left(\log
(1/(1-y))\right)^{\rho_0}}=\sigma_0,\ {\rm for}\ \rho_0>1.
$$

 Note that the point $z=1$ must be an essential
singularity of the function $G$, since the entire function
$F_{r-1}$ interpolating $G$'s coefficients is trascendental.
$\Box$

\section{Proofs of Lemmas \ref{levin2},  \ref{x} and \ref{xx}}

{\bf Proof of Lemma \ref{levin2}} Let  denote
$$\eta=\limsup_{k\to\infty}\frac{\xi(k)}{k}\log|a_k|. $$

First we prove that
$\eta\le(\rho+1)(\sigma_h\rho)^{1/(\rho+1)}/\rho$. If
$\sigma_h=+\infty$, then the inequality is trivial. Suppose
$\sigma<+\infty$ and $\sigma>\sigma_h$. Then $$ \log M(y,h)<\sigma
V(\frac 1{1-y}),\ {\rm for}\ y_0<y<1,$$ which yields
$$\log|a_k|<\sigma V(\frac 1{1-y})+k\log\frac 1y,$$ for $y_0<y<1$
and all $k=0,1,2,\ldots .$

Setting $y=(\xi(k/(\sigma\rho))-1)/\xi(k/(\sigma\rho))$ we obtain
$$\log|a_k|<\sigma V\left(\xi(\frac k{\sigma\rho})\right)+
k\log\left(\frac{\xi(\frac k{\sigma\rho})}{\xi(\frac
k{\sigma\rho})-1}\right)$$ for sufficiently large $k$.

Remembering that $\xi(t)$ is the inverse function of $xV(x)$ we
can write $$\log|a_k|<\frac k{\rho\xi(\frac
k{\sigma\rho})}+\frac{k}{\xi(\frac k{\sigma\rho})}\{1+o(1)\}$$
$$=\frac {\rho+1}\rho\frac k{\xi(\frac k{\sigma\rho})}\{1+o(1)\},\
k\to\infty.$$

By properties (a) and (d) of proximate orders (see
p.~\pageref{prop}) we have that $\xi(t)=\tilde
V(t)=t^{\tilde\rho(t)},\ \tilde\rho(t)\to 1/(\rho+1),\
t\to\infty.$ Hence,
$$\lim_{t\to\infty}\frac{\xi(lt)}{\xi(t)}=l^{1/(\rho+1)}$$ and $$
\xi(\frac
k{\sigma\rho})=\frac{\xi(k)}{(\sigma\rho)^{1/(\rho+1)}}\{1+o(1)\},\
k\to\infty.$$

Therefore, for an arbitrary $\varepsilon>0$ we can assert that
asymptotically $$\log|a_k|<
\frac{\rho+1}\rho(\sigma\rho)^{1/(\rho+1)}\frac
k{\xi(k)}(1+\varepsilon)$$ and, thus,

$$\eta=\limsup_{k\to\infty}\frac{\xi(k)}{k}\log|a_k|
\le\frac{(\rho+1)}\rho(\sigma_h\rho)^{1/(\rho+1)}. $$

Now we want to prove that
$\eta\ge(\rho+1)(\sigma_h\rho)^{1/(\rho+1)}/\rho$. If
$\eta=+\infty$, then the inequality is trivial. Suppose that
$\eta<+\infty$ and $\eta=(\rho+1)(\sigma\rho)^{1/(\rho+1)}/\rho$
for some $\sigma<\sigma_h$. We will obtain a contradiction. Choose
$\sigma_1,\ \sigma<\sigma_1<\sigma_h.$ Then $$ \log|a_k|+k\log
y<k\left\{\frac{\rho+1}\rho\cdot\frac {
(\sigma_1\rho)^{1/(\rho+1)}}{\xi(k)}+\log y\right\},\ k>k_0.$$

In the previous inequality we put $$k=\left[\frac
{\sigma_h\rho}{1-y}V(\frac 1{1-y})\right]$$ meaning $k$ equal to
the entire part of the number between parenthesis and assuming
$1-y$ so small that $k>k_0$. Remembering that $\xi(t)$ is
increasing and a regularly varying function (property (e)), for an
arbitrary $\varepsilon>0$ we can write the following inequality
for $k>k_0$

$$ \log|a_k|+k\log y$$ $$<\frac{\sigma_h\rho}{1-y}V(\frac 1{1-y})
\left(\frac{\rho+1}\rho\left(\frac
{\sigma_1}{\sigma_h}\right)^{1/(\rho+1)}\frac1{\xi(\frac1{1-y}V(\frac
1{1-y}))}-(1-y)\right)(1+\varepsilon)$$ $$=\sigma_h\rho
V(\frac1{1-y})\left(\left(\frac
{\sigma_1}{\sigma_h}\right)^{1/(\rho+1)}\frac{\rho+1}\rho-1\right)(1+\varepsilon),$$
(remember that $\xi(t)$ is the inverse function of $xV(x)$).

Since ${\sigma_1}/{\sigma_h}<1$, we have $$ \log|a_k|+k\log
y<\delta\sigma_hV(\frac1{1-y})(1+\varepsilon),$$ for some
$\delta<1$, which yields $$ \limsup_{y\to
1^-}\frac{\log\mu(y,h)}{V(1/(1-y))}=\sigma_\mu<\sigma_h,$$ where
$\mu(y,h)=\max\{|a_k|y^k:k\in\mathbb N\cup\{0\}\}$.

On the other hand, for any $y'>y$ the inequalities $$ M(y,h)\le
\sum_{k=0}^\infty |a_k|(y')^k\left(\frac y{y'}\right)^k \le\frac
{y'}{y'-y}\mu(y',h)$$ hold. Setting $y'=1-s+sy$, where $s,\
0<s<1,$ is to be chosen later, we obtain
\begin{equation}\label{ineqe} M(y,h)\le\frac
2{(1-s)(1-y)}\mu(1-s+sy,h).\end{equation}

Since $\rho>0$ and $V(x)$ is a regularly varying function it
follows that

$$ \sigma_h=\limsup_{y\to 1^-}\frac{\log M(y,h)}{V(1/(1-y))}\le
\limsup_{y\to 1^-}\frac{\log\mu(1-s+sy,h)}{V(1/(1-y))}$$ $$\le
\limsup_{y\to 1^-}\frac{\log\mu(1-s+sy,h)}{V(1/(1-(1-s+sy)))}
\lim_{y\to 1^-}\frac{V(1/(s(1-y)))}{V(1/(1-y))}$$
$$=\sigma_\mu\left(\frac 1s\right)^\rho.$$

Taking an $s$, such that $(\sigma_\mu/\sigma_h)^{1/\rho}<s<1,$ we
obtain $\sigma_h\le\sigma_\mu,$ which shows that the inequality
$\sigma_\mu<\sigma_h$ obtained earlier, and hence
$\sigma<\sigma_h,$ is impossible. $\Box$

{\bf Proof of Lemma \ref{x} }Let $\omega>\sigma_0,\ \tau>\rho_0,$
then asymptotically
\begin{equation}
\log M(x,F)<\omega(\log x)^\tau. \label{izn}
\end{equation}
Thus, $$ \log|C_n|<\omega(\log x)^\tau-n\log x$$ asymptotically
with respect to $x$ and for $n=0,1,2,\ldots.$

The usual method of finding extrema applied to the right side of
the previous inequality shows that asymptotically
\begin{equation}\label{ineq}
\log|C_n|<\omega(1-\tau)\left(\frac
n{\omega\tau}\right)^{\tau/(\tau -1)}.
\end{equation}

Conversely, assume that the asymptotic inequality (\ref{ineq})
holds. Then $$ |C_nx^n|<K_1\exp\{h(n)\}$$ for $n=0,1,2,\ldots,$
and $x\ge 1,$ where $K_1=K_1(\omega,\tau)$ is a positive constant
and $$h(n)=-Kn^\eta+n\log x,\ \
K=\frac{\omega(\tau-1)}{(\omega\tau)^{\tau/(\tau-1)}}>0,\ \
\eta=\frac \tau{\tau-1}.$$

We next analyze $h(n)$ for real values of its argument. Since
$\eta>1$, the function $h(n)$ attains its maximum equal to
$$K(\eta-1)\left(\frac{\log
x}{K\eta}\right)^{\eta/(\eta-1)}$$ at the point $$\tilde
n=\left(\frac{\log x}{K\eta}\right)^{1/(\eta-1)}.$$

Substituting $K$ and $\eta$ for their expressions in terms of
$\omega$ and $\tau$ we obtain the asymptotic inequality $$ \log
\mu(x,F)\le\omega(\log x)^\tau,$$ where
$\mu(x,F)=\max\{|C_n|x^n:n\in\mathbb N\cup\{0\}\}.$

On the other hand,  for any $x'>x>0$,
$$M(x,F)\le\sum_{n=0}^\infty |C_n|(x')^n\left(\frac x{x'}\right)^n
\le \frac {x'}{x'-x}\mu(x,F).$$ Setting $x'=2x$ we obtain
$$\log M(x,F)\le \omega(\log 2x)^\tau+\log 2$$ and
\begin{equation}\label{mch}\log M(x,F)<
(\omega+\varepsilon)(\log x)^\tau\end{equation} asymptotically
for any arbitrary $\varepsilon
>0$.

Thus (\ref{ineq}) follows from (\ref{izn}) and (\ref{mch}) from
(\ref{ineq}). This shows that the logarithmic order of the
function $F(z)$ is the infimum of the numbers $\tau$ satisfying
(\ref{ineq}) and the logarithmic type of the same function is the
infimum of the numbers $\omega$ satisfying (\ref{ineq}) for $\tau$
equal to $\rho_0$. From this conclusion both assertions of the
theorem follow at once. $\Box$

{\bf Proof of Lemma \ref{xx} }Let $\omega>\sigma_0,\ \tau>\rho_0,$
then asymptotically
\begin{equation}
\log M(y,h)<\omega\left(\log \frac 1{1-y}\right)^\tau.
\label{izn2}
\end{equation}

In the previous inequality, for $k$ large enough, we set
$y=(k-1)/k$, obtaining asymptotically
\begin{equation}\label{ineq2}
\log|a_k|<(\omega+\varepsilon)(\log k)^\tau
\end{equation}
for any arbitrary $\varepsilon>0$.

Conversely, assume that the asymptotic inequality (\ref{ineq2})
holds. Then $$|a_k|y^k<K_1\exp\{h(k)\}$$ for $k=1,2,3,\ldots$ and
$y<1$ close enough to 1, where $K_1$ is a positive constant and $$
h(k)=(\omega+\varepsilon)(\log k)^\tau+k\log y.$$

We next analyze $h(k)$ for real values of its argument $k\ge 1$.
We have $\lim_{k\to\infty} h(k)=-\infty$. Also,
$$h'(k)=\frac {\tau(\omega+\varepsilon)(\log k)^{\tau-1}}k+\log
y,$$ which means that for $y<1$ close enough to 1 we have
$h'(2)>0$. On the other hand, $h'(k)\to \log y<0,\ k\to\infty.$ We
can conclude that, for $y<1$ close enough to 1, the function
$h(k)$ attains its maximum on the interval $[2,\infty)$ at a point
$\tilde k$ such that $h'(\tilde k)=0$. It is easy to prove that
$h'(\tilde k)=0$ yields
$$\tilde k=\frac{\tau(\omega+\varepsilon)}{1-y}\left(\log\frac
1{1-y}\right)^{\tau-1}(1+o(1)) ,\ y\to 1^-,$$ and $$ h(\tilde
k)=(\omega+\varepsilon)\left(\log\frac1{1-y}\right)^\tau(1+o(1)),\
y\to 1^-.$$

Thus,  $$\mu(y,h)<\exp\left\{(\omega+2\varepsilon)\left(\log\frac
1{1-y}\right)^\tau\right\}$$ for any arbitrary $\varepsilon>0$ and
$y$ close enough to $1^-$.

Now using (\ref{ineqe}) and remembering that $\tau>1$ we can write
the following inequalities $$ M(y,h)\le\frac
2{(1-s)(1-y)}\exp\left\{(\omega+2\varepsilon)\left(\log\frac
1{s(1-y)}\right)^\tau\right\}$$ and
\begin{equation}\label{mch2}
M(y,h)\le \exp\left\{(\omega+\varepsilon_1)\left(\log\frac
1{1-y}\right)^\tau\right\}
\end{equation}
for any arbitrary $\varepsilon_1>0$ and $s,\ 0<s<1$ and $y$ close
enough to $1^-$.

Thus, (\ref{ineq2}) follows from (\ref{izn2}) and (\ref{mch2})
from (\ref{ineq2}). This shows that the order of the function
$h(z)$ is the infimum of the numbers $\tau$ satisfying
(\ref{ineq2}) and the type of this  function is the infimum of the
numbers $\omega$ satisfying (\ref{ineq2}) for $\tau$ equal to
$\rho_0$. From this conclusion both assertions of the theorem
follow at once. $\Box$

\vspace{0.5cm}

{\bf Acknowledgements} I would like to express my gratitude to
Prof. I.V. Ostrovskii for reading the manuscript of this paper and
making valuable remarks. I would also like to thank Prof. O.M.
Katkova and Prof. A.M. Vishnyakova for fruitful discussions on
this topic.


\begin{thebibliography}{99}

\bibitem{aissen}
Aissen M., Edrei A., Schoenberg I.J., Whitney A.: On the
Generating Functions of Totally Positive Sequences, {\sl Proc.
Nat. Acad. Sci. U.S.A.,} {\bf 37} (1951), 303--307.

\bibitem{dhpf}
Alzugaray M.T.: Domains of Holomorphy of Generating Functions of
P\'olya Frequency Sequences of Finite Order, accepted for
publication in the journal {\sl Positivity}, Kluwer Academic
Publishers.

\bibitem{beuer}
Beuermann F.:  Wachtumsordnung, Koeffizientenwachstum und
Nullstellendichte bei Potenzreihen mit endlichem Konvergenzkreis,
 {\sl Math. Z.}, {\bf 33} (1931), 98-108.


\bibitem{bieb}
Bieberbach L.: {\sl Analytische Fortsetzung}, Springer-Verlag,
Berlin, 1955.

\bibitem{faber}
Faber G.: Beitrag zur Theorie der ganzen Funktionen, {\sl Math.
Ann.}, {\bf 70} (1911), 48--68.

\bibitem{gelf}
Gelfond A.O.: Sur un th\'eor\`eme de MM. Wigert-Leau, {\sl Matem.
sb.}, {\bf 36} (1929), 99--101.


\bibitem{tp}
Karlin S.: {\sl Total Positivity}, Vol. I, Stanford University
Press, California, 1968.

\bibitem{ind}
Katkova O.M.: {\sl On Indicators of Entire Functions of Finite
Order with Multiply Positive Sequences of Coefficients (Russian)},
Dep. VINITI 22.02.89, No.1179-B89.

\bibitem{growth}
Katkova O.M.: On the growth of entire generating functions of
multiply positive sequences, accepted for publication in the
journal {\sl Matematicheskaya Fizika, Analiz, Geometriya}.


\bibitem{katos}
Katkova O.M., Ostrovskii I.V.: Zero Sets of Entire Generating
Functions of P\'olya Frequency Sequences of Finite Order, {\sl
Math. USSR-Izvestiya,} {\bf 35} (1990), 101--112.

\bibitem{lev}
Levin B.Ya.: {\sl Distribution of Zeros of Entire Functions},
Transl. Math. Monographs, vol.5, AMS, Providence, RI, 1980.

\bibitem{macint}
Macintyre A.J., Wilson R.: Associated Integral Functions and
Singular Points of Power Series, {\sl J. London Math. Soc.,} {\bf
22} (1948), 298--304.

\bibitem{schoenb}
Schoenberg I.J.: On the Zeros of the Generating Functions of
Multiply Positive Sequences and Functions, {\sl Ann. of Math.,}
{\bf 62} (1955), 447--471.

\bibitem{sen}
Seneta E.: {\sl Regularly varying functions}, Springer-Verlag,
Berlin - Heidelberg - New York, 1976.


\bibitem{she1}
\v{S}eremeta M.N.: On the Connection Between the Growth of the
Maximum Modulus of an Entire Function and the Moduli of the
Coefficients of its Power Series Expansion, {\sl Amer. Math. Soc.
Transl.} {\bf 2}, {\bf 88} (1970), 291--301.


\bibitem{she2}
Sheremeta M.N.: O Svyazi mezhdu Rostom Tselykh ili Analiticheskikh
v Krugie Funktsiy Nulevovo Poryadka i Koeffitsientami ikh
Stepennyh Razlozheniy, {\sl Izv. Vyssh. Uchebn. Zav.}, Mat. {\bf
6}(73) (1968), 115--121.

\vspace{0.5cm}

Maria Teresa Alzugaray \\
A.D. de Matem\'atica, Faculdade de Ci\^encias e Tecnologia,
Universidade do Algarve, Gambelas, 8000-117 FARO, PORTUGAL\\
Email: mtrodrig@ualg.pt

\end{thebibliography}
\end{document}